\documentclass{amsart}
\usepackage{amsmath}
\newcommand{\Hilb}{{\rm{Hilb}}}

\newcommand{\codim}{{\rm{codim}}\;}
\newcommand{\coker}{{\rm{coker}}\;}

\newcommand{\Pic}{{\rm Pic}}

\newcommand{\Hom}{{\rm{Hom}}}

\newcommand{\Cat}{{\rm{Cat}}}
\newcommand{\rank}{{\rm{rank}}\;}

\newcommand{\Tor}{{\rm{Tor}}}

\newcommand{\PP}{{\mathbb P}}

\newcommand{\GG}{{\mathbb G}}
\newcommand{\RR}{{\mathbb R}}

\newcommand{\cO}{{\mathcal O}}
\newcommand{\cM}{{\mathcal M}}
\newcommand{\QED}{{\hfill  $\sqcap\!\!\!\!\sqcup$ \medskip}}

\newtheorem{lemma}{\bf Lemma}[section]
\newtheorem{proposition}[lemma]{\bf Proposition}
\newtheorem{theorem}[lemma]{\bf Theorem}
\newtheorem{corollary}[lemma]{\bf Corollary}

\newtheorem{definition}[lemma]{\bf Definition}

\newtheorem{example}[lemma]{\bf Example}

\newtheorem{remark}[lemma]{\bf Remark}
\newtheorem{remarks}[lemma]{\bf Remarks}
\begin{document}

\title[Algebra of Fano 3-folds of genus 12]{Geometry and Algebra of
prime Fano 3-folds of genus 12}
\author[Frank-Olaf Schreyer]{Frank-Olaf Schreyer}

\address{Fakult\"at f\'ur Mathematik und Physik,
Universit\"at Bayreuth,
D-95440 Bayreuth,
Germany}
\email{schreyer@btm8x5.mat.uni-bayreuth.de } 

\date{8.11.1999}
\begin{abstract}
The connection between these Fano 3-folds and plane quartic curves is
explained.
\end{abstract}

\maketitle
\tableofcontents

\section{Introduction}

 According to Mukai \cite{Muk1},\cite{Muk2} every prime Fano 3-fold of 
genus 12 has geometric realizations in three different ways:

\medskip
(a) Let $V$ be a 7-dimensional vectorspace and let  $\eta \colon \Lambda^2 V \to N$  be a net
of alternating forms on $V$, ie. $N$ is 3-dimensional. Denote by 
$$\GG(3,V,\eta) = \{ E \in \GG(3,V) |  \Lambda^2 E \subset \ker(\eta \colon \Lambda^2 V \to N) \}$$ 
the Grassmannian of isotropic 3-spaces of $V$.
\medskip

(b) Let $f \in S_4(U)$ be the equation of a plane quartic $F \subset \PP(U)=\PP^2$. 
A polar hexagon $\Gamma$ of $F$ is the union of six lines 
$ \Gamma = \{ l_1\cdot \ldots \cdot l_6 = 0 \}$ such that
$f = l_1^4 + \ldots +  l_6^4$.
We can also identify $\Gamma$ with a sixtuple of points in $\check\PP^2 = \PP(U^*)$, 
ie. a point of $ \Hilb_6 (\check\PP^2)$. Then the variety of sums of powers presenting f is
$$VSP(F,6) = \overline{ \{ \Gamma \in \Hilb_6(\check\PP^2) | \Gamma
\hbox{ is polar to } F  \} },$$
the closure of the set of polar hexagon of $F$.   
\medskip

(c) Let W be a 4-dimensional vector space, and let $q \colon U^* \hookrightarrow S_2W^*$ be a net of quadrics
in $\check\PP^3=\PP(W^*)$. 
Consider the Hilbert scheme $\Hilb_{3t+1}(\PP^3)=H_1 \cup H_2$ 
and the component $H_1$ containing the twisted cubic curves $C \subset \PP(W) = \PP^3$.
Let
$$H \subset \GG(3,S_2W)$$ 
be the image of $H_1$ under the map

$$\Hilb_{3t+1}(\PP^3) \to \GG(3,H^0(\PP^3,\cO(2))) = \GG(3,S_2W),$$which sends  $$ C \mapsto H^0(\PP^3,I_C(2)).$$
Denote by$$H(q) = H \cap \GG(3,ker(S_2W \to U)) \subset \GG(3,S_2W)$$the variety of twisted cubics, whose quadratic equations are annihilated by $U^* \subset S_2W^*$. 

\begin{theorem}\label{Muk}{\rm (Mukai \cite{Muk1,Muk2}) } Let $X$ be a prime Fano 
3-fold of
genus 12
over an algebraically closed field of characteristic 0. Then
there exist
(a) a net of alternating forms $\eta \colon \Lambda^2 V \to N$,
(b) a plane quartic F,
(c) a net of quadrics $q \colon U^* \hookrightarrow S_2(W)^*$
such that
$$ X \cong \GG(3,V,\eta) \cong VSP(F,6) \cong H(q).$$
Conversely a general net of alternating forms, a general quartic or a general net of
quadrics gives a smooth prime Fano 3-fold of genus 12.
\end{theorem}

\begin{corollary} The moduli spaces
$\cM_{Fano}$ of prime Fano 3-folds of genus 12,
$\cM_3$ of curves of genus 3,$\cM_q$ of nets of quadrics, and $\cM_{3,\vartheta^{ev}}$ of curves of genus 3 together with a non-vanishing
theta characteristic are birational to each other.
\end{corollary}

$\cM_{3,\vartheta^{ev}}$ occurs, since a net of quadrics in $\PP^3$ is determined by its discriminant,
a plane quartic, together with the associated non-vanishing theta characteristic. 
The connection between (a) and (b) is sketched in \cite{Muk1}. 
The surprising fact that $\cM_3$ and $\cM_{3,\vartheta^{ev}}$ are
birational, is actually an old result due to Scorza \cite{Sco1} recently
reconsider by Dolgachev and Kanev \cite{DK}. 
The purpose of this paper is to give a detailed description, how the three realizations are related to
each other.
Inparticular it is proved that Mukai's and Scorza's constructions 
give the same birational transformation. For non-algebraically closed
ground 
fields our investigation gives that the models of type (a) and (b)
exists
over the field of definition of the $V_{22}$, while the space curve
model (c) is in general only defined after a field extension.
\medskip
 
{ \bf Ackowledgement} I am grateful to David Eisenbud, Geir
Ellingsrud, Toni Iarrobino, J\'anos Koll\'ar, Shigeru Mukai, Kristian
Ranestad and
Stein-Arild Str\o mme for fruitful discussions on the material of this paper.

\section{Apolarity and sums of powers}

 Let $k$ be a field of characteristic zero. Consider $S=k[x_0,\ldots,x_r]$ and 
$T=k[\partial_0,\ldots,\partial_r]$. $T$ acts on $S$ by differentiation:
$$\partial^{\alpha}(x^{\beta}) =  \alpha!\frac{\beta}{\alpha} x^{\beta-\alpha}$$
if $\beta \geq \alpha$ and 0 otherwise. Here $\alpha$ and $\beta$ are multi-indices,
${\beta \choose \alpha} = \prod {\beta_i \choose \alpha_i}$ and so on.
Inparticular we have a perfect pairing, {\sl apolarity}, between forms of degree n and 
homogenous differential operators of order n. 

Note that the polar of a form  $f \in S$ in a point $a \in \PP^r$ is given by  $P_a(f)$ for
$a=(a_0,\ldots,a_r) $ and $P_a  = \sum a_i \partial_i  \in T$.

One can interchange the role of S and T by defining
$$x^{\beta}(\partial^{\alpha}) =  \beta!{\alpha \choose \beta } \partial^{\alpha-\beta}.$$
With this notation we have for forms of degree n
$$P_a^n(f)=f(P_a^n)=n!f(a).$$
Moreover $$
f(P_a^m)=0 \iff f(a)=0 \eqno(2.1)$$ if $m \geq n$.

\medskip
Apolarity allows to define artinian Gorenstein graded quotient rings of $T$ 
via forms: For f a homogenous form of degree n define
$$ f^{\bot} = \{D \in T | D(f)=0 \}$$
and $$A^f = T/f^{\bot}.$$
The socle of $A^f$ is in degree n. Indeed 
$P_a(D(f)) = 0 \hskip3pt \forall P_a \in T_1 \iff D(f) = 0 $  or $ D \in T_n$. 
Inparticular the socle of $A^f$ is 1-dimensional, and $A^f$ is indeed Gorenstein.

Conversely for a graded Gorenstein ring $A = T/I$ with socle in degree n
multiplication in A induces a linear form $f\colon S_n(T_1) \to k$ which can be identitied with a homogenous
polynomial $f \in S$ of degree n. This proves:

\begin{theorem}\label{Mac}{\rm (Macaulay, \cite{Mac})}
The map $F \mapsto A^F$ gives a bijection between hypersurfaces $F=\{f=0\} \subset \PP^r$ of 
degree n and Artinian graded Gorenstein quotient rings $A = T/I$ of T
with socle in degree n.
\end{theorem} 

\medskip
Note, that 
$$
(f^{\bot}:D)=D(f)^{\bot}  \eqno(2.2)$$for any homogenous $D \in T$.

\begin{definition}{\rm (cf. \cite{Iar1,Syl3})}
 For forms $f$ of even degree $2n$ the matrix
 
 $$\Cat(f) = (D_iD_j(f))_{1 \leq i,j \leq {n+r \choose r}}$$
 with $D_1,\ldots,D_{n+r \choose r} \in T_n$ a basis, is called catalectican matrix of $f$.
 $f$ is called non-degenerate, if $Cat(f)$ has maximal rank.
\end{definition}

\medskip
The Hilbert function and syzygies of $A^F$ depends on suttle 
 properties of $F$, cf. \cite{Iar1,IK}. 
 For example for plane quartics, we have (cf. Clebsch \cite{Cle}):
 
\begin{theorem}\label{rankThm} 
 Let $F=\{f=0\}$ be a plane quartic. The following integers are equal:
\item[(1)] $\dim_k A^F_2$ ,
\item[(2)] $\rank \Cat(f)$,
\item[(3)] the minimal s such that over the algebraic closure $\overline{k}$ of $k$ $f$ lies in the 
 closure of forms $ l_1^4+\ldots+l_s^4$ .
\end{theorem}
 
{\it Proof.}  $\dim_k A^f_2 = \rank \Cat(f)$ holds because multiplication in $A^F$ gives
 perfect pairings $A^F_2 \times A^F_2 \to A^F_4 \cong k$. 
 
 Suppose  $f = l_1^4+\ldots+l_s^4$. The corresponding lines  $L_1,\ldots,L_s \in \check\PP^2$, 
 viewed as points in the dual space, impose at most s conditions on quadrics 
 $\{D=0\} \subset  \check\PP^2$. Hence $\dim_k f^{\bot}_2 \geq 6-s$ and $\dim_k A^F_2 \leq s$.
 Since the Hilbert function of $A^F$ varies semi-continously with $F$, $\dim_k A^F_2 \leq s$ holds
 for all forms in the closure of the set of forms of the sums of s powers.
 
 Conversely, suppose $\dim_k A^F_2 = s$. $A^F=T/I$ is Gorenstein of codimension 3. Hence the 
 structure theorem of Buchsbaum-Eisenbud \cite{BE} applies: $A^F$ has syzygies
 $$  0 \longleftarrow A^F \longleftarrow T \longleftarrow
 F_1 \buildrel \phi \over \longleftarrow F_2
 \longleftarrow T(-7) \longleftarrow 0, $$
 with $F_1 = \bigoplus_{i=1}^{2r+1} T(-a_i), F_2^* \cong F_1(7),\phi$ skew-symmetric, and $I$
 generated by the $2r \times 2r$ pfaffians of $\phi$. Conversely, any sufficiently general
 skew-symmetric homomorphism $\phi \in \Hom_T(F_2,F_1)$ defines via its pfaffians
 a graded Artinian Gorenstein ring $A=T/I$ with the same Hilbert function as $A^f$. Therefor
 it suffices to establish the sum presentation $f = l_1^4+ \ldots + l_s^4$ for an $f$
 corresponding to an $A=A^F$ with sufficiently general syzygy matrix $\phi$ for each possible
 numerical type of syzygies.
 
 There are only a few number of numerical cases: 
 $$ \begin{matrix} \hbox{\rm Hilbert function} & (a_1,\dots,a_{2r+1}) & m & p
& (b_{r+1},\dots,b_{2r+1})\cr
 (1,3,6,3,1)& (3,3,3,3,3,3,3) &3&0&(4,4,4)\cr
 (1,3,5,3,1) & (2,3,3,3,3) &2&0&(4,4)\cr
 (1,3,4,3,1) & (2,2,3) &0&1&(4)\cr
 (1,3,4,3,1) & (2,2,3,3,4) &1&1&(3,4)\cr  
 (1,3,3,3,1) & (2,2,2,4,4) &0&2&(3,3)\cr 
 (1,2,3,2,1) & (1,3,3) &1&0&(4)\cr
 (1,2,2,2,1) & (1,2,4)&0&1&(3) \cr
 (1,1,1,1,1) & (1,1,5) &0&1&(2) 
\end{matrix}$$ 
 We argue in each of the cases seperately but similarly: Let $n$ be the number of cubic generators
of $I=f^\bot$ and $m=\lfloor \frac{n}{2} \rfloor$.  Consider the $n \times n$ submatrix 
$\tilde\phi$ of $\phi$ corresponding to the linear coefficients of the quartic syzygies.
 
Suppose there is a  $m \times m$ (skew) symmetric submatrix of zeroes in $\tilde\phi$. 
The corresponding quartic syzygies and all the syzygies which involve only equations of degree 
$\leq 2$ give $r=m+p$
syzygies of degrees $(b_{r+1},\dots,b_{2r+1})$ as 
indicated above between equations of degrees $(a_1,\dots,a_{r+1})$.
Here $p$ is the number
of syzgies which involve only equations of degree 
$\leq 2$, and  $m+p=r$. Thus we obtain a block 
decomposition
$$ \phi = \begin{pmatrix} \begin{matrix}0&* \cr -*&0 \end{matrix}
 &\psi \cr  -\psi^t&0 \cr \end{pmatrix}
  \begin{matrix}
 \}  & r+1 \cr \}  & r \cr 
\end{matrix}$$
with a $(r+1) \times r$ matrix $\psi$. The $r+1$ minors of $\psi$ are among the pfaffians of $\phi$;
they are precisely the generators involved in our $r$ syzygies. 
By Hilbert-Burch \cite[Thm 20.15]{Eis}  these minors generate the homogenous ideal $J_{\Gamma}$ of a set
  $\Gamma \subset \check\PP^2$ of distinct points with syzygies
$$  0 \longleftarrow J_{\Gamma} \longleftarrow   
\bigoplus_{i=1}^{r+1} T(-a_i)
 \buildrel \psi \over \longleftarrow \bigoplus_{j=r+1}^{2r+1} T(-b_j) \longleftarrow 0,$$
if $\psi$ is sufficiently general. $J_{\Gamma}$ is 2-regular, since $b_j \leq 4$ for $j \ge r+1$.
 Hence
the Hilbert function of $R=R_{\Gamma} = T / J_{\Gamma}$ takes the values
$$h_R(t)= \dim_k (R_{\Gamma})_t = \deg \Gamma$$ for $t \geq 2$. On the other hand
$$ \dim_k (R_{\Gamma})_2 = \dim_k A^f_2=s$$
as $(J_{\Gamma})_{\leq 2} = I_{\leq 2}$ by construction. Thus for sufficiently general $\psi$ 
$\Gamma$ consist of s points
$L_1=\{l_1=0\},\ldots,L_s=\{l_s=0\}$. To prove that there exists a sum presentation 
$f=\lambda_1l_1^4+\ldots+\lambda_sl_s^4$ we consider $T \to R \to A$ and the induces inclusions
$$\Hom(A_4,k) \subset \Hom(R_4,k) \subset \Hom(T_4,k).$$ The linear forms
$$ \{D \mapsto D(l_i^4) \}$$
are contained in $\Hom(R_4,k)$. Moreover since $\Gamma$ imposes $s$ independent conditions
on quartics, these linear forms span the image. Inparticular $$\{D \mapsto D(f) \} \in \Hom(A_4,k)$$
is contained in this space, ie.
$D(f) = D(\lambda_1l_1^4+\ldots+\lambda_sl_s^4)$ for all $D \in T_4$
for suitable $\lambda_1,\ldots,\lambda_s \in k$. Hence $f=\lambda_1l_1^4+\ldots+\lambda_sl_s^4$
as desired. Taking roots of the $\lambda_i$'s we can put them into the equations
$l_i$. 

It remains to prove the  existence of a $m \times m$ block of zeroes in $\tilde\phi$,
 possibly after row and column operations.
 Let $V_f:= \Tor_2^T(A^f,k)_4$. Then $\tilde\phi$
corresponds to a net of alternating forms $\Lambda^2V_f \to T_1$ and we are looking
for a subspace $E \in \GG(m,n)=\GG(m,V_f)$ such that 
$\Lambda^2E \subset ker(\Lambda^2V_f \to T_1)$. If $m=1$ there is nothing to prove. If $m\geq2$
then $E$ excists, because for $j={m \choose 2}$ and $c_j=c_j(\Lambda^2 {\mathcal E}^*)$ the $j^{th}$
 Chern class, where ${\mathcal E}$ denotes the universal subbundle on $\GG(m,n)$,
we have $c_j^3 \neq 0$.
\QED 

\medskip
\noindent
Notice that we expect a 3-dimensional family in case $s=6$, a 1-dimensional family of sum
 presentations for $f$ with $s=5$, or $s=4$ and $r=3$, or $s=3$ and $h_A(1)=2$,  and a unique
presentation otherwise.

\begin{remarks}{\rm 1) The fact that, despite the dimension count, a general plane 
quartic is not the sum of 5 powers, goes back to Clebsch \cite{Cle}.

\noindent
 2) It is not true, that every $f$ with $\dim_k A^f_2=s$ is a sum of $s$ powers. For $s$ with an
 unique presentation examples are rather obvious. But even in case $s=5$ this occurs: Eg.
 $$f=(1-\frac{1}{t^2})x_1^4+x_1^3(x_0-\frac{4}{t}x_2)+\frac{1}{t^2}(x_1+tx_2)^4
+x_2^3(x_0-4tx_1)+(1-t^2)x_2^4$$
is not the sum of 5 powers. The reason is that the quadric in $I=f^\bot$ is a double line. Hence
distinct points in a $\Gamma$ would give a linear form in $I$, a contradiction.  For $f$
the 1-dimensional family of sum decompositions degenerates to the family parametrized by $t$ 
of decompositions into five summands as above.}
\end{remarks}

\begin{definition}
For $F=\{f=0\} \subset \PP^n$ a hypersurface we call a scheme $X
\subset \check \PP^n$ apolar to $F$ if $I_X \subset F^\bot$. The
family
of zero-dimensional apolar subschems of degree $s$ of $F$ is denoted by
$VPS(F,s)$. 
\end{definition}
Note that with this definition
$$VSP(F,s) \subset VPS(F,s)$$
is an open subscheme and equality holds if $VPS(F,s)$ is irreducible
and
$VSP(F,s)$  non-empty.

\begin{theorem}\label{ab} Let $F=\{f=0\} \subset \PP(U)$ be a non-degenerate plane quartic.
Then $$VPS(F,6) \cong \GG(3,V_f,\eta_f),$$ where
$V_f = (f^{\bot})_3^*$,$N_f=U^*$ and $\eta \colon \Lambda^2 V_f \to N_f$ the skew-symmetric syzygy matrix of
$A^f$. Conversely, for a net $\eta \colon \Lambda^2 V \to N$ of skew-forms on a 7-dimensional
vector space, whose pfaffians define a ideal $I$ of codimension 3 in $S(N)$, the dual socle 
quartic $F=F(V,\eta) \subset \PP(N^*)$ is a non-degenerate quartic, and 
$VPS(F(V,\eta),6) \cong \GG(3,V,\eta).$ 
\end{theorem}

{\it Proof.} By \ref{Mac}, the structure theorem of Buchsbaum-Eisenbud
\cite{BE}, and  \ref{rankThm}
$$ F \mapsto (\Lambda^2 V_f \to \eta_f)$$
and $$(\eta \colon \Lambda^2 V \to N)  \mapsto F(V,\eta)$$ give bijections between
$$\{F | \det(\Cat(f)) \ne 0 \} \longleftrightarrow \{\eta \colon \Lambda^2 V \to N | \codim I = 3 \}.$$
Moreover points $p \in \GG(3,V,\eta)$ correspond to block decompositions 
$$ \phi =   \begin{pmatrix} \begin{matrix}0&*\cr-*&0 \end{matrix} &\psi \cr
  -\psi^t&0 \end{pmatrix}  \begin{matrix}\} & 4 \cr \} & 3 \end{matrix}$$
of the syzygy matrix $\phi$.
We claim that for every $\psi$ the ideal of minors
$I(\psi)$ has codimension 2, hence defines of a subscheme of length 6 in $\check\PP^2$. 

Assume that, $I(\psi)$ has not depth 2. Then by Hilbert-Burch, the corresponding minors 
have a common factor. Since the minors are minimal generators of $I_{pf}$, the factor has to be a 
linear form $t \in T_1$. So $\psi$ is a matrix of syzygies among 4 quadrics without a common factor. The
quadrics generate an ideal $J$ of codimension $\ge 2$. Let $B=T/J$. $B$
has Hilbert function $(1,3,2,1,\ldots)$. If $ \dim B = 0$ then 3 general quadrics in J
form a regular sequence, whose quotient has Hilbert function $(1,3,3,1,0)$ and the fourth
quadric cuts down to a ring with Hilbert function $(1,3,2,0)$. This is not the case.  
 So $\dim B =1$ and $B$ has Hilbert function $(1,3,2,1,1,\ldots)$. 
Such quotients B exist: B is defined by 4 quadrics in the homogenous ideal of a point $p \in \check\PP^2$.
 However such a $\psi$ does not occur as part of a skew-symmetric
matrix $\phi$, whose pfaffians have codimension 3.
The syzygies of B start
$$  0 \longleftarrow B \longleftarrow T \longleftarrow
 4T(-2)  \longleftarrow 3T(-3) \oplus T(-4)\oplus \ldots \longleftarrow \ldots . $$
 Since $tJ \subset I_{pf}$ the syzygy
$4T(-2) \leftarrow T(-4)$ gives a relation among the pfaffians. But this relation is not
in the space generated by the columns of $\phi$, since the sequence
$$ 3T(-3) \buildrel -\psi^t \over \longleftarrow 4T(-4) \longleftarrow T(-6) \longleftarrow 0$$
is exact. This contradicts the exatness of the pfaffian complex.

Thus we have a well-defined morphism
$\alpha\colon \GG(3,V,\eta) \to \Hilb_6(\check\PP^2).$
To prove that $\alpha$ is
 an isomorphism
onto its image, consider the open part $\Hilb_6(\check\PP^2)^o$ of the Hilbert scheme of length 6
subschemes, which impose independent conditions on quadrics, and the embedding
$\Hilb_6(\check\PP^2)^o \hookrightarrow \GG(4,T_3)$. The diagram
$$\begin{matrix}\Hilb_6(\check\PP^2)^o & \hookrightarrow &\GG(4,T_3) \cr
 \alpha \uparrow & & \uparrow \cr
\GG(3,V,\eta) & \hookrightarrow & \GG(4,V^*) \cr
\end{matrix}$$
commutes, where $V^* = (f^{\bot})_3 \subset T_3$.

Finally, note that the image of $\alpha$ contains all polar hexagons of $F$. Indeed, if
$f = l_1^4+ \ldots +l_6^4$ for distinct lines $\Gamma = \{L_1,\ldots,L_6\} \subset \check \PP^2$, then
$\Gamma$ imposes independent conditions on quadrics by Thm \ref{rankThm} Hence syzygies of $\Gamma$ are
of type:
$$ 0 \longleftarrow R_{\Gamma} \longleftarrow T \longleftarrow 4T(-3) \buildrel \psi \over 
\longleftarrow 3T(-4) \longleftarrow 0. \eqno(2.3)$$
By (2.1) the ideal $J_{\Gamma} \subset f^{\bot}$. Hence we have a sequence
$$ 0 \longleftarrow A^f  \longleftarrow R_{\Gamma} \longleftarrow I_{A/R} 
\longleftarrow 0. \eqno(2.4)$$
Since $A^f$ and $R=R_{\Gamma}$ have Hilbert functions $(1,3,6,3,1)$ and $(1,3,6,6,6,\ldots)$
respectively, $I_{A/R}$ has 3 cubic generators with 4 linear relations:
$$0 \longleftarrow I_{A/R} \longleftarrow 3T(-3) \longleftarrow 4T(-4).$$
The minors of the presentation matrix are contained in the annihilator, which is $J_{\Gamma}$.
Hence, this matrix is $\psi^t$ again, and $I_{A/R} \cong \omega_R(-4)$. A mapping cone between
the complex (2.3) and its dual over the sequence (2.4), gives syzygies of $A^f$ with the
desired block structure.

We do not claim at this point, that every non-degenerate $f$ has a non-degenerate polar hexagon.
However, if there is one, then the points in the image of $\alpha$ corresponding to them,
form an open subset.  
\QED

\begin{corollary} For a general plane quartic $F$ the variety of polar 
hexagons $VSP(F,6)$ is a smooth Fano 3-fold of genus 12.
\end{corollary} 

{\it Proof.} Since for the tautological subbundle ${\mathcal E}$ on $\GG(3,V)$ the sheaf $ \Lambda^2 {\mathcal E}^*$ is 
globally generated by $ \Lambda^2 V^*$, a general net $N$ of skew-forms defines a smooth
 subscheme $\GG(3,V,\eta)$ of codimension 9. Since $\omega_{\GG(3,V)} \cong {\cO}_{\GG(3,V)}(-7)$ and
$\Lambda^9(3\Lambda^2{\mathcal E}^*) \cong {\cO}_{\GG(3,V)}(-6)$ one has 
$\omega_{\GG(3,V,\eta)} \cong {\cO}_{\GG(3,V,\eta)}(-1)$.
By degree reasoning  $\GG(3,V,\eta)$ is irreducible and hence it
 is a Fano 3-fold.  
\QED

\section{The Scorza map}

 In this section we recall some results of Scorza from \cite{DK}.
\medskip

A plane cubic $C$ is called anharmonic, if $C$ lies in the PGL(3)-orbit closure of 
$\{x_0^3+x_1^3+x_2^3=0\}$. The reason is that the cross-ratio of the Fermat cubic is anharmonic.
Let ${\bf A} \subset \PP(H^0(\PP^2,\cO(3))^*) \cong \PP^9$ denote the variety of anharmonic cubics. The 
PGL-orbit closure of a general cubics is a hypersurface of degree 12. 
Due to the additional automorphism of the Fermat cubic, ${\bf A}$ is hypersurface of degree 4. 
In terms of coordinates $(a,\dots,j)$ of $\PP^9$, 
 $ ax_0^3+bx_1^3+cx_2^2+3dx_0^2x_1+3ex_0^2x_2+3fx_1^2x_0+3gx_1^2x_2+3hx_2^2x_0+
3ix_2^2x_1+6jx_0x_1x_2,$
${\bf A}$ is defined by the Aronhold invariant
$$\begin{matrix} {\bf I}_4 = & abcj-(bcde+cafg+abhi)-j(agi+bhe+cdf) \cr
 & +(afi^2+ahg^2+bdh^2+bie^2+cgd^2+cef^2) \cr
 & -j^4+2j^2(fh+id+eg)-3j(dgh+efi) \cr
 & -(f^2h^2+i^2d^2+e^2g^2) +(ideg+egfh+fhid).  
\end{matrix}$$

\begin{lemma}\label{anhcub}  Let $g \in k[x_0,x_1,x_2]$ be a plane cubic. The following are equivalent:
\vskip0pt (1) $A^g \cong T/(D_1,D_2,D_3) $ is a complete intersection of three quadrics,
\vskip0pt (2) g is not an anharmonic cubic.
\end{lemma}

{ \it Proof.}  If $g$ is not a cone, then $A^g$ has Hilbert function $(1,3,3,1)$, and there
are precisely three quadrics in $I=g^{\bot}$. 
 By \cite{BE} either 
$A^g$ is a complete intersection of three quadrics, or $A^f$ has syzygies
$$  0 \longleftarrow A^g \longleftarrow  T \longleftarrow  
\bigoplus_{i=1}^{5} T(-a_i)
 \buildrel \psi \over \longleftarrow \bigoplus_{j=1}^{5} T(-b_j) \longleftarrow T(-6)
\longleftarrow 0,$$
with $(a_1,\ldots,a_5) = (2,2,2,3,3)$ and $b_i = 6-a_i$. In the second case the three
quadrics of $I=f^{\bot}$ intersect in 3 points  $\{L_1,L_2,L_3\} \in
\check\PP^2$ (possibly infinitesimal near), and as in section 2 we obtain
$g=\lambda_1l_1^3+\lambda_2l_2^3+\lambda_3l_3^3$. Conversely, if $g$ is a smooth
 anharmonic cubic all three quadrics of $I$ vanish in  
$\{L_1,L_2,L_3\} \in \check\PP^2$, and $A^g$ is not a complete intersection. 

Since all cones are anharmonic cubics, this proves the lemma.  
\QED

Let $F= \{f=0\} \subset \PP^2$ be a non-degenerate plane quartic.  Consider
$$S_{F} = \{a \in \PP^2| P_a(F) \in {\bf A} \}.$$
Then either $S_{F} = \PP^2$ or $S_{F}$ is a plane quartic. The first case does not occur for
non-degenerate quartics, \cite[6.6.3]{DK}. We call $S_{F}$ the covariant quartic of $F$.
Consider
$$T_{F} = \{ (a,b) \in \PP^2 \times \PP^2 | \rank P_{a,b}(F) \leq
1\}.$$
\begin{lemma}{\rm (\cite[6.8.1]{DK})}
 Let $F \subset \PP^2$ be a general quartic. Then $S_{F}$ is a smooth quartic and 
$T_{F}$ is a smooth symmetric correspondence of type $(3,3)$ on $S_{F} \times S_{F}$ 
without united points. 
\end{lemma}

{  \it Proof.} For a complete proof see \cite{DK}. The  reason, why $T_{F}$ is such a correspondence
on $S_{F} \times S_{F}$ is the following:

 Suppose $(a,b) \in T_{F}$, say
$P_{a,b}(f) = h^2$. Set $t_0=P_b \in T_1$ and $(t_1,t_2) = (h^{\bot})_1 \subset T_1$. Then
$t_0t_1,t_0t_2 \in (P_a(f)^{\bot})_2$ and $(P_a(f)^{\bot})$ is not a complete intersection of
quadrics. By Lemma \ref{anhcub} $P_a(F)$ is an anharmonic cubic. So $a \in S_{F}$. 

For general $F$ and general $a \in S_{F}$, $P_a(F)$ is a smooth Fermat cubic, and then the points
 $b \in \PP^2$ such that $\rank P_{a,b}(F) \leq 1$ are the 3 vertices of the Hessian triangle of
$P_a(F)$. So $T_{F}$ is symmetric of type (3,3). 
Since 
$$ 2{\partial^2 f \over \partial x_i \partial x_j}(a) ={\partial^2 \over \partial x_i \partial x_j} P_{a,a}(f)$$
by (2.1), $T_{F}$ has no united points, if
 $\rank( {\partial^2 f \over \partial x_i \partial x_j}(a)) \geq 2$ for all $a \in S_{F}$. 
This is the case  for general $F$, because then the Hessian 
$He(F) = \{\det({\partial^2 f \over \partial x_i \partial x_j})=0 \}$
 is smooth. 
\QED

  Let $S$ be a smooth plane quartic and $\vartheta$ an even theta characteristic
on $S$. Consider the $\vartheta$-correspondence
$$T_{\vartheta} = \{ (a,b) \in S \times S | h^0(S,\vartheta(a-b)) \geq 1\}.$$
$\vartheta$ is not effective. 
 $deg$ $ \vartheta(a) = 3$, and 
$h^0(S,\vartheta(a-b)) = h^1(S,\vartheta(a-b)) = h^0(S,\vartheta(b-a)).$ 
 So $T_{\vartheta}$ is a symmetric correspondence of type (3,3) without united points.

\begin{theorem}{\rm (\cite[7.6]{DK} )} Let k be a algebraically closed field.
 If $F$ is a quartic such that $S_{F}$ is a smooth quartic, then there exists
a unique theta characteristic $\vartheta = \vartheta_{F}$ on $S_{F}$ such that
$$T_{F} = T_{\vartheta}.$$ 
\end{theorem}

{ \it Proof.} For a general $F$ and a general $(a,b) \in T_{F}$ consider the polar Hessian
triangles 
$T_{F}(a) = b+b_1+b_2$ and $T_{F}(b)= a+a_1+a_2$, ie. the vertices of the Hessian of 
$P_a(F)$ and $P_b(F)$ respectively. All 6 points are different. If $P_{a,b}(f)=h^2$ then both
$b_1,b_2$ and $a_1,a_2$ span the line $\{h=0\}$. So $S_{F} \cap \{h=0\} = \{a_1,a_2,b_1,b_2\}$
and
$$T_{F}(a)-a+T_{F}(b)-b=b_1+b_2+a_1+a_2 \eqno(3.1)$$
is a canonical divisor on $S_{F}$. Moreover 
$$T_{F}(a)-a \equiv T_{F}(b)-b \eqno(3.2)$$ for any 2 points $a,b \in S_{F}$.
To establish this consider the map
$$ S_{F} \to \Pic^2(S_{F}), a \mapsto \cO(T_{F}(a)-a). $$
Since even for degenerate polar Hessians the three (possibly infinitesimal near) points
 $T_{F}(a)=b_1+b_2+b_3$ are not colinear, 
$h^0(S_{F},\cO(T_{F}(a))) = 1$ for all $a \in S_{F}$. So $$h^0(S_{F},\cO(T_{F}(a)-a)) = 0,$$ as
$T_{F}$ has no united points. It follows, that the image of $S_{F}$ does not intersect the 
$\Theta$-divisor of $\Pic^2(S_{F})$. Since $\Theta$ is ample, $S_{F}$ maps to a point. 
By (3.1),(3.2) $$\vartheta = \cO(T_{F}(a)-a))$$
is an non-vanishing theta characteristic.
\QED

\begin{remark}{\rm Although the $[\vartheta] \in \Pic \, C$ is a point defined
over the ground field $k$, the line bundle $\vartheta$ may not be
defined
over the ground field. For an example see \ref{realexample} below.}
\end{remark}

\begin{theorem} {\rm (Scorza, \cite[7.8,7.11]{DK})} Let k be algebraically closed.
 The rational map induced by
$$ s\colon F \mapsto (S_{F},\vartheta_{F})$$
from the moduli spaces of curves of genus 3
$$ s\colon \cM_3 \to  \cM_{3,\vartheta^{ev}}$$
to the moduli of curves of genus 3 together with a even theta
characteristic is birational.
\end{theorem}

\begin{remarks} {\rm
(1) The projection $ \cM_{3,\vartheta^{ev}} \to \cM_3 $ is a finite cover of degree
36:1, since a curve of genus 3 has precisely 36 even theta characteristics.

(2) For a pair $(S,\vartheta)$ of a plane quartic together with a non-vanishing 
theta characteristic, the quartic $F=s^{-1}(S,\vartheta)$ is called Scorza quartic of
 $(S,\vartheta)$.  Dolgachev and Kanev give two desription of 
$s^{-1}$. In section 5 we will give another one.

(3) More general, for a canonical curve $S \subset \PP^{g-1}$ of genus $g$, and a non-vanishing 
theta characteristic $\vartheta$ with some plausible, but yet unproven hypothesis
Scorza constructs a quartic hypersurface $F \subset \PP^{g-1}$. See \cite{Sco2} and \cite{DK}.}
\end{remarks}

\section{Rank 2 vector bundles on $\PP^3$ with $c_1=0$ and $c_2=3$}

  Let $W$ be a 4-dimensional vector space, $\check \PP^3=\PP(W^*)$ and 
$q \colon U^* \hookrightarrow S_2W^*$ a net of quadrics in $\check \PP^3$, whose general 
element is smooth. Let
$$W \to U \otimes W^*$$
the associated symmetric matrix with entries in U, and
$$ b_q \colon W \otimes \cO_{\PP^2}(-2) \to W^* \otimes \cO_{\PP^2}(-1)$$
the associated map of sheaves on $\PP^2 = \PP(U)$ twisted. Let 
$S = S_q = \{\det(b_q )=0 \} \subset \PP^2$ 
denote the discriminant of the net, and $\vartheta = \vartheta_q = \coker (b_q)$.	Since
$b_q$ is symmetric,
$$\vartheta \cong {\mathcal E}xt^1_{\cO_{\PP^2}}(\vartheta,\omega_{\PP^2}) \cong {\mathcal H}om_{\cO_S}(\vartheta,\omega_S).\eqno(4.1)$$		           
If $S$ is smooth, then $\vartheta$ is an invertible $\cO_S$-module, hence $\vartheta$ is a non-vanishing
theta characteristic on $S$. 

Conversely, given a plane quartic $S$ and a torsion free rank 1 $\cO_S$-module $\vartheta$ which
satisfies (4.1), we denote by $W^* = H^0(S,\vartheta(1))$ and 
$$q = q(S,\vartheta) \colon U^* \hookrightarrow S_2W^*$$
the corresponding net of quadrics in $\PP(W^*) = \check \PP^3$.

If S is smooth, then 
$$ \phi_{\vartheta(1)} \colon S \hookrightarrow \check \PP^3$$
is an embedding, and the image
$ \tilde S=\phi_{\vartheta(1)}(S)$ is the variety of vertices of the singular quadrics in the net.
Equation of $ \tilde S \subset\check \PP^3 = \PP(W^*)$ are given by the $3 \times 3$ minors of
$$\tilde b_q \colon W \otimes \cO_{\check\PP^3}(-1) \longrightarrow U \otimes \cO_{\check\PP^3}.$$
If $S$ is not smooth, one can take these equations to define $\tilde S$.
From another point of view, $\tilde b_q$ is the jacobian matrix of q. 

\bigskip We consider now the apolarity pairing between $\PP^3 = \PP(W)$ and 
$\check \PP^3 = \PP(W^*)$. By $R=SW$ and $T=SW^*$ denote the homogenous coordinate rings 
respectively. Let
$$ q^{\bot} = \{D \in R \hskip1pt | D(Q)=0 \hskip1pt \forall Q \in q(U^*) \subset S_2W^*\}, $$ and
$$ A^q = R/q^{\bot}.$$ 
$A^q$ is an Artinian ring with Hilbert function $(1,4,3,0, \ldots)$. $A^q_1 = W$, $ A^q_2 = U$ and
multiplication given by $q\colon S_2W \to U$. 

\begin{lemma}\label{syzAq}  Let $q$ be a general net of quadrics. Then 
$A^q$ has syzygies:
{\small $$  0 \leftarrow A^q \leftarrow R \leftarrow
 7R(-2) \buildrel (\phi_1,\phi_2) \over \longleftarrow 8R(-3) \oplus 3R(-4)
 \leftarrow 8R(-5) \leftarrow 3R(-6) \leftarrow 0. $$ }
\end{lemma}

 {\it Proof.} The  number of syzygy $\dim \Tor^R_i(A^q,k)_j$ in the above
sequence take the minimal possible values for an Artinian ring $A$ with
Hilbert function $(1,4,3,0,\ldots)$. Thus by semi-continuity it suffices to establish
 the existence of one example $q$ with such syzygies. The Kleinian net
$$q_{\scriptstyle Klein}=({\scriptstyle{1\over2}}z_1^2-z_0z_2,{\scriptstyle{1\over2}}z_2^2-z_0z_3,{\scriptstyle{1\over2}}z_3^2-z_0z_1)$$
where $T = k[z_0,z_1,z_2,z_3]$, has this property.
\QED

 Consider the map $\phi_1$ in the complex above, and its kernel sheafivied and 
twisted by $\otimes \cO_{\PP^3}(5) $:
$${\mathcal E}={\mathcal E}_q= ker(7\cO_{\PP^3}(3) 
\buildrel \phi_1 \over \longleftarrow 8\cO_{\PP^3}(2)).$$

\begin{proposition}\label{syzE} Let $q$ be a general net of quadrics. Then  ${\mathcal E}_q$
is a stable rank 2 vectorbundle with Chern classes $c_1=0,c_2=3$ and syzygies
$$ 0 \longleftarrow {\mathcal E}_q \longleftarrow 8\cO_{\PP^3}(-2) \longleftarrow
 7\cO_{\PP^3}(-3)  \longleftarrow \cO_{\PP^3}(-5)  \longleftarrow 0. \eqno(4.2)$$ 
Its $H^2$-cohomology module is
$$ A^q(5) = \bigoplus_{n} H^2(\PP^3,{\mathcal E}_q(n)).$$
\end{proposition}

{\it Proof.} Since $A^q$ Artinian the first syzygy module sheafivied is a rank 6 vector bundle,
$$  0 \longleftarrow \cO_{\PP^3}(5) \longleftarrow
 7\cO_{\PP^3}(3)  \longleftarrow {\mathcal F} \longleftarrow  0.\eqno(4.3)$$
${\mathcal E}_q$ is the kernel of by $\phi_1$,
$$ (0 \longleftarrow) {\mathcal F} \longleftarrow  8\cO_{\PP^3}(2) \longleftarrow {\mathcal E}_q \longleftarrow  0. \eqno(4.4)$$
One expects that $ {\mathcal F} \longleftarrow  8\cO_{\PP^3}(2)$ is surjective outside a set of
codimension $8-6+1=3$, hence that ${\mathcal E}_q$ is a vector bundle of rank 2 outside a finite
set of points. The expected number of these points is 0 by Porteous formula.
Thus either $\phi_1$ has rank $\leq 5$ along at least a curve, or ${\mathcal E}_q$ is a  rank 2 
vector bundle with Chern polynomial
$$c_t({\mathcal E}_q) = {(1+5t)(1+2t)^7 \over (1+3t)^8} \equiv 1+3t^2 \qquad mod \quad t^4.$$ 
For a general $q$ the second alternative takes place, as one can check by considering 
an example, eg. the Kleinian net.

Since ${\mathcal E}_q$ has rank 2 and $c_1=0$, wedge product 
$${\mathcal E}_q \otimes {\mathcal E}_q \to \Lambda^2 {\mathcal E}_q \cong \cO_{\PP^3}$$ 
gives ${\mathcal E}_q^* \cong {\mathcal E}_q$. Thus the dual of the sequences (4.3) and (4.4) give the
exact sequence (4.2). Since this complex is short enough to stay exact on global sections
for arbritrary twists, this is the minimal resolution. The last statement follows from the 
cohomology sequence of (4.2.2) and (4.2.3). ${\mathcal E}$ is stable, because $H^0(\PP^3,{\mathcal E})=0$ 
, cf. \cite[Lemma II 1.2.5]{OSS}.
\QED

\begin{corollary}\label{dualities}  If $q$ is general, then there are natural isomorhism
$$ U \cong \Tor^R_4(A^q,k)_6 \cong (\Tor^R_2(A^q,k)_4)^*$$
and a skew-symmetric self-duality
$$ \Tor^R_3(A^q,k)_5 \cong (\Tor_3^R(A^q,k)_5)^*.$$
The maps $3R(-4) \leftarrow 8R(-5)$ and $8R(-5) \leftarrow 3R(-6) $ in the resolution
(4.2) are dual to each other under these isomorphisms. 
\end{corollary}

{\it Proof.} The matrices yield minimal presentations
$$0 \leftarrow \bigoplus_n H^1(\PP^3,{\mathcal E}(n)) \leftarrow 3R(1) \leftarrow 8R$$
and
$$0 \leftarrow Ext^4_R(A^q,R) \leftarrow 3R(6) \leftarrow 8R(5)$$
respectively. By Serre duality and $A^q(5) =\bigoplus_n H^2(\PP^3,{\mathcal E}(n))$,
$$Ext^4_R(A^q(5),R)= \bigoplus_n H^1(\PP^3,{\mathcal E}^*(n)).$$
Since ${\mathcal E} \cong {\mathcal E}^*$, the modules are isomorphic, hence the desired isomorphisms
follow by comparison of the presentations. The self-duality on $\Tor^R_3(A^q,k)_5$ is skew,
since the isomorphism ${\mathcal E} \cong {\mathcal E}^*$ has this property. Finally, we note
$$\begin{matrix}U  &\cong  A^q_2 \cong H^2(\PP^3,{\mathcal E}(-3)) \cong H^1(\PP^3,{\mathcal E}^*(-1))^* \cr
  & \cong H^1(\PP^3,{\mathcal E}(-1))^* \cong \Tor^R_2(A^q,k)_4^* \cong
\Tor^R_4(A^q,k)_6.
\end{matrix}  $$
\QED

\begin{corollary} 
 $\tilde S=\tilde S_q \subset \check \PP^3$ is the variety of unstable planes of ${\mathcal E}_q$.  
$\tilde S$ determines ${\mathcal E}_q$ up to isomorphism. The moduli space $\cM_{\PP^3}(2;0,3)$ of
 rank 2 vector bundles on $\PP^3$ with $c_1=0$ and $c_2=3$ has a component birational to
$\GG(3,S_2W^*)$.
\end{corollary}

{\it Proof.} A plane $H = \{h=0\} \subset \PP^3$ is unstable for ${\mathcal E}$, iff
$H^0(H,{\mathcal E}\mid_H) \ne 0$, equivalently, if multiplication with $h$ is not injective on
$$ H^1(\PP^3,{\mathcal E}(-1)) \buildrel h \over \longrightarrow H^1(\PP^3,{\mathcal E}).$$ 
$H^1(\PP^3,{\mathcal E}(-1)) \cong U^*$ and $H^1(\PP^3,{\mathcal E}) \cong W^*$. A quadric
$q_1 \in U^* \subset S_2W^*$ is annihilated by $h \in W$, iff 
$Q_1=\{q_1=0\} \subset \check \PP^3$ is a cone with vertex $H \in  \check \PP^3$. So the variety
of unstable planes coincides with the variety $\tilde S \subset\check \PP^3$ of vertices of the 
cones. $\tilde S$ determines $q$, which in turn determines $A^q$ and ${\mathcal E} _q$. The vector
bundles obtained from points $q \in \GG(3,S_2W^*)$ form an open part of the moduli scheme
$\cM_{\PP^3}(2;0,3)$, since by semi-continuity and minmality, the cohomology modules 
$\bigoplus_n H^2(\PP^3,{\mathcal E}(n)) $ have the same numerical type of syzygies for an open part
of $\cM_{\PP^3}(2;0,3)$.
\QED

\section{Twisted cubics annihilated by a net of quadrics}

 Let $q \colon U^* \hookrightarrow S_2W^*$ be a net of quadrics as before. 
Let $H(q)$ denote the variety of twisted cubics $C \subset \PP^3$, whose equations 
$H^0(\PP^3,{\mathcal I}_C(2)) \subset S_2W$ are annihilated by $q$. Let 
$V_q = (q^{\bot})_2 \subset S_2W$. Since a twisted cubic is defined by its quadrics and
$h^0(\PP^3,{\mathcal I}_C(2)) = 3$,
$H(q)$ is a subset of $\GG(3,V_q)$ in a natural way. We are looking for an net of alternating 
forms on $V_q$, which defines $H(q) \subset \GG(3,V_q)$. 

For a description of $\Hilb_{3t+1}(\PP^3)$  and the map
$$\Hilb_{3t+1}(\PP^3) \to \GG(3,S_2W)$$ see \cite{EPS,PS}. 

\medskip
\noindent Consider the syzygies of $A^q$. By definition of $A^q$ we have
$$\Tor^R_1(A^q,k)_2 \cong V_q.$$
Define
$$N_q = \Tor^R_2(A^q,k)_4 $$
and consider 
$$\eta_q \colon \Lambda^2 V_q \rightarrow N_q$$
given by multiplication in the algebra $\Tor^R_*(A^q,k)$,
cf. \cite{Eis}[Exercise A3.20] 
$$ \eta(p_1\wedge p_2) = 0$$
 for $p_1,p_2 \in V_q$, iff the Koszul syzygy 
$$p_1 \otimes p_2-p_2 \otimes p_1 \in ker(R \leftarrow 7R(-2))$$ 
lies in $Im(7R(-2) \leftarrow 8R(-3))$.

\begin{theorem} \label{Hq=GGVq} 
$$H(q) \cong \GG(3,V_q,\eta_q)$$
for a general net $q$.
\end{theorem}

{\it Proof.} Let $C \subset \PP^3$ be a rational normal curve
whose ideal $I_C=(p_1,p_2,p_3) $ is generated by three quadrics $p_1,p_2,p_3 \in V$.
$C$ has syzygies
$$  0 \longleftarrow \cO_C \longleftarrow \cO_{\PP^3} \longleftarrow
 3\cO_{\PP^3}(-2)  \longleftarrow 2\cO_{\PP^3}(-3) \longleftarrow  0.\eqno(5.1)$$
Hence all syzygies among $p_1,p_2,p_3$ are generated by linear relations, and
$E = (I_C)_2 = (p_1,p_2,p_3)_2 \in \GG(3,V,\eta)$. 

Conversely, suppose that
$E=(p_1,p_2,p_3) \in \GG(3,V,\eta).$  We will prove that $p_1,p_2,p_3$ generate 
the homogenous ideal of a curve $C$ of degree 3 and arithmetic genus 0. 
Choose $p_4,\ldots,p_7 \in V$, such that $p_1,\ldots,p_7$ form a basis. By definition of $\eta$
there is a matrix $8R(-3) {\buildrel \psi \over \leftarrow} 3R(-4)$, such that
$$\phi_1 \cdot  \psi = 
\begin{pmatrix} 0 & -p_3 & p_2 \cr
p_3 & 0 & -p_1 \cr
-p_2 & p_1 & 0 \cr
0 & 0 & 0 \cr
\vdots && \vdots \cr \end{pmatrix}
$$
gives the matrix of Koszul relations, where $\phi_1$ is the matrix of
linear syzygies in Lemma \ref{syzAq}.
With $$\tau =  \begin{pmatrix} p_1 \cr p_2 \cr p_3 \end{pmatrix}$$
we have 
$$\phi_1 \cdot \psi \cdot \tau = 0.$$
But  for $ ker(7\cO_{\PP^3}(-2)  \buildrel \phi_1 \over \longleftarrow 8\cO_{\PP^3}(-3)) = {\mathcal E}(-5)$
we have $H^0(\PP^3,{\mathcal E}(1))=0$ by Prop. \ref{syzE}. So
$$ \psi \cdot \tau = 0,$$
ie. $p_1,p_2,p_3$ are three quadrics with some linear relations
$$r_1 \cdot p_1 + r_2 \cdot p_2 + r_3 \cdot p_3 = 0.\eqno(5.2)$$
The coefficients $r_1,r_2,r_3 \in R_1$ are linearly independent for a general linear 
combination of rows of $\psi$. Because otherwise any two elements of $E$ would have a 
common linear factor, which implies, that all three elements $p_1,p_2,p_3$ have a common factor,
and  $\Tor^R_3(A^q,k)_4 \ne 0$. But this group is zero by Lemma \ref{syzAq}. Thus 
$$(p_1,p_2,p_3) = \Lambda^2 \tau_2 $$
for $2 \times 3$ matrix 
$$\tau_2 = \begin{pmatrix} r_1 & r_2 & r_3 \cr r_4 & r_5 & r_6 \end{pmatrix}$$
of linear forms $r_i \in R_1$, (cf. \cite{Sch}, Lemma 4.3). Since the minors $p_1,p_2,p_3$
have no common factor, the Hilbert-Burch complex of $\tau_2$ is exact, and $p_1,p_2,p_3$
generate the ideal of a curve $C \subset \PP^3$ of degree 3 and arithmetic genus 0. 
$C \in \Hilb_{3t+1}(\PP^3)$ lies in the component $H_1$, which contains the twisted cubics, 
(cf. \cite{PS,EPS}).

Note, that  boundary points corresponding to plane nodal cubics with an embedded 
point at the node, do not occur, since all $C$ are arithmetically
Cohen-Macaulay.
\QED

\begin{proposition} \label{Fq} For a general net $q \colon U^* \hookrightarrow S_2W^*$ of 
quadrics the pfaffians of the net of alternating forms $\eta_q \colon \Lambda^2 V_q \to N_q$ 
defines an Artinian Gorenstein ring with Hilbert function $(1,3,6,3,1)$ with a smooth
dual socle quartic $F_q = F(V_q,\eta_q)$.
\end{proposition}

{\it Proof.} Since the desired property is an open condition on nets $q$ of quadrics,
it suffices to exhibit an example. For the Kleinian net $q_{Klein}$ we obtain as dual socle quartic
$F_{Klein} = \{x_0^3x_1+x_1^3x_2+x_2^3x_0 = 0\} $. $F_{Klein}$ is smooth for $char(k) \ne 7$.
For  $char(k) = 7$ one can take some other example.
\QED

\section{The Hilbert schemes of lines on $X$}
 
 Let $F=\{f=0\} \subset \PP^2=\PP(U)$ be a non-degenerate plane quartic.
In this section we prove that the circle of constructions

 $$\begin{matrix} F \mapsto (S_F, \vartheta_F),\qquad  &\hbox{ Scorza, } \cr
(S,\vartheta) \mapsto q_ {S,\vartheta },\qquad  &\hbox{  net corresponding to $\vartheta$,}\cr
 q \mapsto A^q,\qquad  &\hbox{ apolarity,} \cr
 A^q \mapsto (\eta_q \colon  \Lambda^2 V_q \to N_q),\qquad  &\hbox{ Tor multiplication,}\cr
 (\eta \colon  \Lambda^2 V \to N) \mapsto A_ {V,\eta },\qquad  &\hbox{ pfaffians,}\cr
 A \mapsto F_A ,\qquad  &\hbox{ dual socle quartic,}
\end{matrix}$$
gives the identity transformation on an open set of quartics. Note that, since
$N_q \cong U^*$ by Cor. \ref{dualities}, $F_A$ is again a quartic in $\PP(U)$. 

 For $X=\GG(3,V,\eta)$ we denote by
${\mathcal H}_X$ the Hilbert scheme of lines  
in $X$ with respect to the Pl\"ucker embedding  $X \subset \GG(3,V) \hookrightarrow 
 \PP(\Lambda^3 V^*).$
 
\begin{theorem} Let q be a general net of quadrics in $\check \PP^3$. Let
${\mathcal E}_q$ be the corresponding vector bundle on $\PP^3$, $X=H(q)=\GG(3,V_q,\eta_q)$ and $F=F_q$
the dual socle quartic of $A_{V_q,\eta_q}$. The following curves are isomorphic:
\vskip0pt (a) the discriminant $S_q$ of $q$,
\vskip0pt (b) the variety $\tilde S_q$ of unstable planes of ${\mathcal E}_q$,
\vskip0pt (c) the Hilbert scheme ${\mathcal H}_X$ of lines on $X$,
\vskip0pt (d) the covariant quartic $S_F$ of $F$. 
\end{theorem}

{\it Proof.}$ (a) \leftrightarrow (b):$ $S_q \cong \tilde S_q$ holds
 by section 4.
 $(b) \leftrightarrow (c) : $ Let $H=\{r=0\} \subset \PP^3$ be an unstable plane. Then 
$$U^* = H^1(\PP^3,{\mathcal E}_q(-1)) \buildrel r \over \rightarrow  W^*=H^1(\PP^3,{\mathcal E}_q)$$ 
is not injective, equivalently,
$$\mu_r \colon W \buildrel r \over \rightarrow U$$ 
not surjective. So $ker(\mu_r)$ is at least 2 dimensional, ie. there are 2 elements $r_1,r_2 \in R_1$
such that  $p_1 = r  \cdot r_1,p_2 = r \cdot r_2 \in V_q $. 

$(r_1,r_2) \cap V_q \subset S_2W$ is at least $7-3=4$ dimensional. So there are further 2 quadrics
$$p_3 = a_1  r_1 + a_2  r_2, p_4 = b_1  r_1 + b_2  r_2 \in (r_1,r_2) \cap V_q.$$
Let $C_{(\alpha:\beta)} \subset \PP^3$ be the curve (!) defined by
$$(p_1,p_2,\alpha p_3+\beta p_4) = \Lambda^2 \begin{pmatrix} r_2 & -r_1 & 0\cr
 \alpha a_1+\beta b_1 & \alpha a_2+\beta b_2 & -r \end{pmatrix}.$$
By Thm \ref{Hq=GGVq}
$$ p_1 \wedge p_2 \wedge (\alpha p_3+\beta p_4) \in \GG(3,V_q,\eta_q) \subset 
\PP(\Lambda^3 V_q^*),\hskip5pt (\alpha:\beta) \in \PP^1,\eqno(6.1)$$
 gives a point in ${\mathcal H}_X$. 

Conversely, every line in $H(q)$ is of type (6.1) for some $p_1,p_2,p_3,p_4 \in V_q$. 
$p_1$ and $p_2$ have a common factor $r$, since 
$$\bigcup_{(\alpha:\beta)} C_{(\alpha:\beta)} \subset \{p_1=p_2=0\} \subset \PP^3,$$
and $H=\{r=0\}$ is a unstable plane.

$(c) \rightarrow (d):$ $F =\{f=0\}$ is non-degenerate by Cor. \ref{Fq} and
Thm. \ref{rankThm}. Moreover
$$\GG(3,V_q,\eta_q) = \GG(3,V_f,\eta_f) \cong VSP(F,6).$$
From this point of view a line 
$$ p_1 \wedge p_2 \wedge (\alpha p_3+\beta p_4) \in \GG(3,V_f,\eta_f) $$
corresponds to a syzygy matrix
$$\phi = \begin{pmatrix} 0 & 0 & 0 & 0 & a_{15} & a_{16} & a_{17} \cr
					                            0 & 0 & 0 & 0 & a_{25} & a_{26} & a_{27} \cr
                                      0 & 0 & 0 & a_{34} & a_{35} & a_{36} & a_{37} \cr
0 & 0 & -a_{34} & 0 & a_{45} & a_{46} & a_{47} \cr
-a_{15} & -a_{25} &-a_{35} & -a_{45} & 0 & a_{56} & a_{57} \cr
-a_{16} & -a_{26} &-a_{36} & -a_{46} & -a_{56} & 0 & a_{67} \cr
-a_{17} & -a_{27} &-a_{37} & -a_{47} & -a_{57} & -a_{67} & 0 \end{pmatrix},\eqno(6.2)$$
and the family of submatrices
$$ \psi_{(\alpha:\beta)} = \begin{pmatrix} 0 & a_{15} & a_{16} & a_{17} \cr
   0 & a_{25} & a_{26} & a_{27} \cr
 a_{34} & \alpha a_{35}+\beta a_{45} & \alpha a_{36}+\beta a_{46} &
\alpha a_{37}+\beta a_{47} \end{pmatrix}$$
corresponds to 1-parameter family of polar hexagons with 3 fixed lines $\in \check\PP^2$ defined by
$$  \begin{pmatrix}   a_{15} & a_{16} & a_{17} \cr
    a_{25} & a_{26} & a_{27}  \end{pmatrix}$$
and three moving lines through the common point $a_{34} \in \PP^2$. Hence
the polar $P_{a_{34}}(F) \in {\bf A}$, ie. $a_{34} \in S_F \subset \PP^2$.  

$(c) \leftarrow (d):$ Conversely given $a \in S_F$. Since $F$ is non-degenerate $g=P_a(f)$ is not 
a cone. By Lemma \ref{anhcub} $A^g$ is not a complete intersection. Hence there are three
quadrics $b_1,b_2,b_3 \in (g^{\bot})_2$ with precisely 2 linear syzygies. Then
$ab_1,ab_2,ab_3 \in (f^{\bot})_3$, and the 2 linear syzygies give 2 of the columns of $\phi$
with many zeroes. Thus this gives a decomposition of $\phi$ of shape (6.2). We only have to
check that none of the 6 possibly non-zero entries of the 2 columns can lie on the diagonal.
Suppose 1 or 2 entries lie on the diagonal. Then, since
$b_1,b_2,b_3$ are the minors of the $3 \times 2$ matrix, either one quadrics is zero, or
they have a common factor and a further syzygy. Both cases are impossible.
Thus every point $a \in S_F$ gives a uniquely determined decomposition of 
$\phi$ of shape (6.2). Hence a well-defined point in ${\mathcal H}_X$. 
Since $a=a_{34}$ in this correspondence this is the inverse of $(c) \rightarrow (d)$.

Notice, that under the isomorphisms of section 5 and \ref{dualities} 
$$\eta_q(p_3 \wedge p_4) = a = a_{34} \in N \cong U^* \cong H^1(\PP^3,{\mathcal E}(-1)).$$
Moreover $(p_1,p_2,p_3,p_4)$ have the relations
$$ \tilde \phi_1 =\begin{pmatrix} -r_2 & a_1 & b_1 & 0 \cr
					    r_1 & a_2 &b_2 & 0 \cr
					    0 & -r & 0 & -p_4 \cr
					    0 & 0 & -r & +p_3  
\end{pmatrix}.$$
 $\tilde \phi_1 \cdot \rho = 0$ for 
$$ \rho =\begin{pmatrix} a_1b_2-b_1a_2 \cr p_4 \cr -p_3 \cr  r \end{pmatrix}.$$
Thus, if we choose a basis for $8R(-3) \oplus 3R(-4)$ in the complex
of Cor. \ref{syzAq} with these four 
relation corresponding basis elements, $\rho$ gives one column of the matrix
$$8R(-3) \oplus 3R(-4) \longleftarrow 8R(-5).$$
Since 
$$0 \longleftarrow \bigoplus_n H^1(\PP^3,{\mathcal E}(n)) \longleftarrow (3R(-4) \longleftarrow 8R(-5)) \otimes R(5)$$ 
is a presentation, we obtain, that $\eta_q(p_3 \wedge p_4) = a \in H^1(\PP^3,{\mathcal E}(-1))$ 
is annihilated by $r$. Since $h^0(H,{\mathcal E}|_H)=1$ for $H=\{r=0\}$ any unstable plane,
there is only one quadric cone $Q \in U^*$ with vertex $H \in \check \PP^3$ up to scalars.
Thus under the isomorphisms of Cor. \ref{dualities}, section 5 and $(a) \leftrightarrow (d)$ the curves
$$S_q, S_F \in \PP(U)$$
are actually equal.
\QED

Let $X=X_q = \GG(3,V_q,N_q)$. Denote by
$$T_{\mathcal H} = \{(L_1,L_2) \in {\mathcal H}_X \times {\mathcal H}_X | L_1 \cap L_2 \ne \emptyset, L_1 \ne L_2 \}$$
the correspondence of intersecting lines in $X$.

%, by
%$B_{\mathcal H} \rightarrow {\mathcal H}$ the universal family of lines, and  by 
%$B_X = \bigcup_{L \in {\mathcal H}} L \subset X$ its image in $X$.
%$$\begin{matrix}B_{\mathcal H} & \rightarrow  &B_X \subset X \cr
%\downarrow & & \cr
% {\mathcal H}_X & & \end{matrix}. $$

\begin{corollary}  Let $q$ be a general net of quadrics. Then $B_X$ consists of 
all singuar twisted cubics in $H(q)$, and $S_F \subset \PP(U)$
is the set of the tripel points of polar hexagons to $F=F_q$.
The correspondences
 \vskip0pt (a) $\qquad T_{\vartheta_q} \hbox{ on } S_q ,$ 
 \vskip0pt (b) $\qquad T_{\mathcal H} \hbox{ on } {\mathcal H}_X, $
\vskip0pt (c) $\qquad  T_F \hbox{ on } F_q $ \hfill
are isomorphic. 
\end{corollary}

{\it Proof.} From the proof of the theorem we see, that the curves 
$C_{(\alpha:\beta)} \in H(q)$ on a line 
$$ \{p_1 \wedge p_2 \wedge (\alpha p_3+\beta p_4)\}_ {(\alpha:\beta) \in \PP^1} \in {\mathcal H}\eqno(6.1)$$
are all singular, and that they have the component $\{r_1=r_2=0\} \subset \PP^3$ in common.
Conversely, if a curve $C \in H(q)$ is singular, it is reducible and one of its components
is a line $\{r_1=r_2=0\}$ in the intersection of two reducible quadrics $p_1=r_1 \cdot r,
p_2=r_2 \cdot r \in H^0(\PP^3,I_C(2)) \subset V_q.$ $r$ defines a unstable plane of ${\mathcal E}_q$
and gives a line (6.1). 

Now take the point of view from polar hexagons. If a point $\Gamma = \{L_1,\ldots,L_6 \}$ 
lies on a line in $VSP(F,6)$,
then three lines of the hexagon pass through a common point $a_{34} \in \PP^2$. Conversely, if
$\{L_1,L_2,L_3\} \subset \Gamma$ pass through a point $a \in \PP^2$, then $P_a(F)$ is anharmonic,
ie. $a \in S_F$ and $\tilde f = f - (\lambda_4 l_4^4+\lambda_5 l_5^4+\lambda_6 l_6^4)$
is a quartic with $h_{A^{\tilde f}}(1) = 2$ and $h_{A^{\tilde f}}(2) = 3$, since $f$ is not a sum of
five powers. So by Thm \ref{rankThm} there is a pencil of 3-tuples of lines presenting $\tilde f$, and this
gives the family $\Gamma_{(\alpha:\beta)}$ defined by (6.1). For three values
of $(\alpha:\beta) \in \PP^1$ one of the moving lines passes through an intersection point $b$ of a pair
of the fixed lines. This corresponds to an intersection of two different lines in ${\mathcal H}$, ie.
a point in $T_{\mathcal H}$, and also to the point $(a,b) \in T_F$.

Finally to prove $T_F \cong T_{\vartheta_q}$ note that we already know $S_F = S_q$. Thus
$T_F$ and $T_{\vartheta_q}$ correspond both to one of the 36 even theta characteristics on
$S_F$.  Since $\cM_3$ and $\cM_{3,\vartheta^{ev}}$ have the same dimension, this implies, that 
the circle of constructions induces a covering transformation of
$$\cM_{3,\vartheta^{ev}} \rightarrow \cM_3$$
over an open set corresponding to general nets $q$. Since $\cM_q$ hence also $\cM_{3,\vartheta^{ev}}$
is irreducible, it suffices to varify $T_F = T_{\vartheta_q}$ in one example, where all
steps are defined. For example one can check this for the Kleinian net 
$q_{Klein}$. 

Note, that $F_{q_{Klein}}=S_{q_{Klein}}$ is the Klein curve. Thus,
$\vartheta_{Klein}$ is the unique theta characteristic on the Klein curve invariant
under the whole automorphism group $G_{168}$, cf. \cite[\S 232]{Bur}.
\QED

\begin{corollary} Over an algebraically closed field the circle of 
constructions
 $\quad F \mapsto (S_F, \vartheta_F),\quad 
(S,\vartheta) \mapsto q_ {S,\vartheta },\quad 
 q \mapsto A^q,\quad  
 A^q \mapsto (\eta_q \colon  \Lambda^2 V_q \to N_q),\quad 
 (\eta \colon  \Lambda^2 V \to N) \mapsto A_ {V,\eta },\quad 
 A \mapsto F_A   $
gives the identity transformation on an open set of quartics.
\end{corollary}

\begin{corollary}  Over an algebraically closed field the circle of constructions  define birational transformations
of the moduli spaces $\cM_{Fano}$ of nets of alternating forms (equivalently of prime Fano 3-folds
of  genus 12 by Mukai's Theorem),
$\cM_3$ of curves of genus 3, $\cM_{3,\vartheta^{ev}}$ of curves of genus 3 together with a non-vanishing
theta characteristic, and $\cM_q$ of nets of quadrics.
\end{corollary}

\begin{corollary} Let K be an arbitrary field of characteristic 0 and
$X$ a smooth prime Fano 3-fold of genus 12. Then the Grassmanian model $\GG(3,V,\eta)$ and the plane model
$VSP(F,6)$ are defined over K. 
\end{corollary}
{\it Proof.} By Mukai's Theorems \cite{Muk1,Muk2} all 3 models are
defined over the algebraic closure of K. However the Hilbert scheme
${\mathcal H}_X$ of
lines on $X$ is defined over K, and so is the correspondence $T_{\mathcal H}$ 
of intersecting lines. Identifying $T_{\mathcal H}=T_F$ and $S_F = {\mathcal H}$ in its canonical embedding, we obtain a quadratic
system of equations for the coefficients of the defining equation
$F =\{f=0\}$ with coefficients
in K. So the  by Mukai's result unique solution with a non-degenerate
quartic is defined over K. The equivalence of the Grassmannian model and
space model is defined over K. So also $\eta$ is defined over K.
\QED

The space model $H(q)$ is in general not defined over the ground
 field, due
 to the descent from $T_\vartheta$ to $\vartheta$.
For the real numbers we note:

\begin{remark}{\rm 1) Let $X \cong VSP(F,6)$ be a prime Fano 3-fold of
genus 12 defined over $k$ with smooth covariant quartic $S_F$. 
If the covariant quartic $S_F$ contains a $k$-rational point, then
the space model $H(q)$ is defined over $k$. Indeed if $a \in S_F$ is
defined over $k$,
 then the fiber $T_F(a)$ and line bundle 
$\vartheta = \cO(T_F(a)-a)$ are defined over $k$.

2) The quadrics $$q_0=w_0^2+w_1^2-w_2^2-w_3^2, q_1=w_0w_2+w_1w_3,$$
$$q_2=(w_1+w_2+w_3)^2+(w_0+w_1-w_3)^2-(w_0+w_1+w_3)^2-(w_0+w_1-w_2)^2$$
span a net $q$, whose discriminant $S_q$ contains no real point.
Indeed the catalectican $\Cat(S_q)$ is positive definit. So the
existence of a point is sufficient but not necessary for the existence
of the space model, even for the ground field $\RR$.
}
\end{remark}

\begin{example}\label{realexample} {\rm For the
Mukai-Umemura
quartic cf. \cite{MU} there are 2 different plane models over $\RR$,
$$F_{MU}=\{(x^2+y^2+z^2)^2=0\} \hbox{ and } 
F^\prime_{MU} = \{ (x^2+y^2-z^2)^2=0 \}.$$
Their covariant quartics are equal to themselves.

For the space  model $H(q_{MU})$ there is only one version. The net of quadrics $q_{MU}$ is the ideal
of the twisted cubic, which, if defined over $\RR$, is isomorphic to 
$\PP^1_\RR$. Its discriminant is the indefinit $F^\prime_{MU}$.
Thus for $VSP(F_{MU},6)$ there is no $\RR$-isomorhic space model

Since
the image of $\{q=(q_0,q_1,q_2) \hbox{ over } \RR\} \to \{(\eta \colon \Lambda^2
V \to \RR^3) \hbox{ over } \RR\} \to \{F \subset \PP^2 \hbox{ over }
\RR \}$  is closed (in a neighborhood
of $F_{MU}$) and $F_{MU}$ is not in the image, we obtain that
 for any quartic $F$ over $\RR$ nearby $F_{MU}$,  the Fano 3-fold
$VSP(F,6)$ has no $\RR$-isomorphic space model. 
}
\end{example}

\end{document}